\documentclass{article}
\usepackage[english]{babel}
\usepackage[utf8]{inputenc}
\usepackage{amsmath}
\usepackage{graphicx}
\usepackage{amsfonts}
\usepackage{enumitem}
\usepackage{amssymb}
\usepackage[colorinlistoftodos]{todonotes}
\usepackage{geometry}
\usepackage{amsthm}
\usepackage{enumitem} 
\usepackage{hyperref}

\newcommand{\remove}[1]{}
\usepackage{tikz}

\long\def\onefigure#1#2{
\begin{figure*}[tbp]
\begin{center}
#1
\end{center}
\caption{#2}
\end{figure*}
} 

\newcommand{\lipefig}[2]  
{\onefigure{\mbox{\psfig{file=#1.eps}}}{\label{f:#1} #2} }

\newcommand{\norm}[1]{\left\lVert#1\right\rVert}

\newcommand{\floor}[1]{\left\lfloor #1 \right\rfloor}
\newcommand{\ceil}[1]{\left\lceil #1 \right\rceil}
\newcommand{\paren}[1]{\left( #1 \right)}

\newcommand{\set}[1]{\left\{ #1 \right\}}
\newcommand{\setcond}[2]{\left\{ #1 : #2 \right\}}

\newcommand{\one}{\mathbf{1}}

\newcommand{\eps}{\varepsilon}

\newcommand{\RR}{\mathbb{R}}

\newcommand{\ZZ}{\mathbb{Z}}

\newtheorem{thm}{Theorem}[section]
\newtheorem{lem}[thm]{Lemma}

\newtheorem{claim}[thm]{Claim}
\newtheorem{fact}[thm]{Fact}

\theoremstyle{definition}
\newtheorem{defn}[thm]{Definition}

\DeclareMathOperator{\conv}{conv}
\DeclareMathOperator{\cl}{cl}

\title{Balancing games on unbounded sets}
\author{Imre Bárány \thanks{Alfr\'ed R\'enyi Institute of Mathematics,13-15 Re\'altanoda Street, Budapest, 1053 Hungary, and Department of Mathematics, University College London Gower Street, London WC1E 6BT. Email: {\tt imbarany@gmail.com}. Research partially supported by Hungarian National Research Grants No 131529, 131696, and 133819.}
\and Jeck Lim\thanks{Mathematical Institute, University of Oxford. Email: 
{\tt jeck.lim@maths.ox.ac.uk}. Research supported by ERC Advanced Grant 883810. Both authors are partially supported by the NSF under Grant No. DMS-1928930, while the authors were in residence at the Simons Laufer Mathematical Sciences Institute in Berkeley, California, during the Spring 2025 semester.}}
\date{}

\begin{document}
\maketitle

\begin{abstract} For a finite set $V\subset \RR^n$, a set $T\subset \RR^n$ is called $V$-closed if $t \in T$ and $v\in V$ imply that either $t+v\in T$ or $t-v \in T$. The set $P(V):=\{\sum_{v \in W} v: W \subset V\}$ is clearly $V$-closed and so are its translates. We show, assuming $V$ contains no parallel vectors, that if $T$ is closed and $V$-closed, and $x \in T$ is an extreme point of $\cl \conv T$, then there is a translate of $P(V)$ containing $x$ and contained in $\conv T$. This result is used to determine the value of a special balancing game. A byproduct is that when $m\ge 2$ and is not a power of 2, then the $m$-sets of a $2m$-set can be coloured Red and Blue so that complementary $m$-sets have distinct colours and every point of the $2m$-set is contained in the same number of Red and Blue sets.
\end{abstract}

\section{Introduction}

A balancing game is a perfect information game played by two players; the Pusher (he) and the Chooser (she). Given a set $V \subset \RR^n$ and a convex set $K\subseteq \RR^n$, the balancing game $G(V,K)$ proceeds through a sequence of points $z_0=0,z_1,z_2, \ldots \in \RR^n$. In round $k$ of the game Pusher offers a vector $v_k \in V$ and Chooser selects $\eps_k \in \{-1,1\}$ and the next position is $z_{k+1}=z_k+\eps_kv_k$. Pusher wins if he can make $z_k\not\in K$ for some $k$ and Chooser wins if she can maintain $z_k\in K$ for all $k$. Balancing games in various forms were introduced by Spencer~\cite{Spen} in 1979. 

\smallskip
In this paper, we consider the case when $V \subset \RR^n$ is finite. 
\begin{defn}
A set $T \subset \RR^n$ is called \emph{$V$-closed} if $t \in T$ and $v\in V$ imply that either $t+v\in T$ or $t-v \in T$.  
\end{defn}

\medskip
Define the set $P(V):=\{\sum_{v \in W} v: W\subset V\}$. $P(V)$ is clearly $V$-closed and so are its translates. Moreover, $P(V)$ is centrally symmetric with respect to its center which is the same as its center of gravity. We denote it by $g(V)=\frac 12 \sum_{v\in V}v$. In a balancing game one often has to check if $g(V) \in P(V)$ or not. According to \cite{Karp} (see also \cite{Garey}) it is NP-hard to decide if $g(V) \in P(V)$ even if $V$ consists of positive integers.

\smallskip
The connection between $V$-closed sets and balancing games is based on the following fact. If some $z_i$ lies in a $V$-closed set $T$, then Chooser can guarantee that $z_k\in T$ for all $k\ge i$. This is basically the only strategy that Chooser can use to win. 

\begin{thm}[{\cite[Theorems 1 and 2]{Bar79}}]\label{th:determined}
    The game $G(V,K)$ is determined:
    \begin{enumerate}
        \item Pusher has a winning strategy if and only if there is no $V$-closed set $T$ with $0\in T\subseteq K$.
        \item Chooser has a winning strategy if and only if there is a $V$-closed set $T$ with $0\in T\subseteq K$.
    \end{enumerate}
\end{thm}

We derive some properties of $V$-closed sets under the condition that $n\ge 2$ and
\begin{equation}\label{eq:nonpar}
 V \subset \RR^n \mbox{ is finite and contains no parallel vectors}.
\end{equation}
In particular we show in Section~\ref{sec:proof} that, under certain natural conditions, a translate of $P(V)$ is contained in a $V$-closed set. Recall that for a convex set $C\subset \RR^n$, a point $x\in C$ is an \emph{extreme point} if there is no distinct pair of points $x_0,x_1\in C$ such that $x$ lies on the open line segment connecting $x_0$ and $x_1$. A simple fact (proved by induction on $n$) from convex geometry is that a nonempty closed convex set in $\RR^n$ contains no line iff it contains an extreme point. We remark that if $T$ is $V$-closed, then its closure, $\cl T$ is also $V$-closed. This explains the condition ``$T$ is closed" in the following theorem. 

\begin{thm}\label{th:vclos} Assume $V$ satisfies condition (\ref{eq:nonpar}), $T\subset \RR^n$ is closed and $V$-closed. Let $C=\cl \conv T$ be the closure of $\conv T$. Assume $x$ is an extreme point of $C$. Then there is a translate $t+P(V)$ of $P(V)$ with $x\in t+P(V) \subset C$.
\end{thm}

\begin{figure}[h!]
\centering
\includegraphics[scale=1.1]{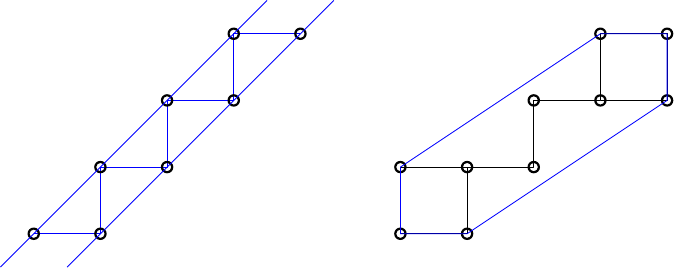}
\caption{Examples of $V$-closed sets.}
\label{fig:line}
\end{figure}

As $x$ is an extreme point of, $C$ contains no line which is important for the application of this theorem to balancing games. When $T$ is bounded, $C$ can't contain a line and in that case Theorem~\ref{th:vclos} is a slightly stronger form of a lemma in \cite[page 120]{Bar79}. The condition that $x$ is an extreme point of $C$ is also necessary as shown by the example on Figure~\ref{fig:line} 
left, where $V=\{e_1,e_2\}.$ (Here and later $e_i,\; i\in [n]$ form the standard basis of $\RR^n$.) Then $T$ is the infinite set of lattice points in the strip $\{(x_1,x_2\in \RR^2: 0\le x_2-x_1 \le 1\}$, and $T$ contains no translate of $P(V)$. In the example on the the same figure, right, $V$ is the same set and the set $T$ of the circled points is $V$-closed and some points of $T$ are not contained in any translate of $P(V)$ that is contained in $C$.

\medskip
The following unusual balancing game $G(V,K_M)$ is considered in the USA January Team Selection Test for IMO 2017~\cite{artof}. The set of vectors Pusher can choose from is $V=\{-1,1\}^n$ and $K_M=\{(x_1,\ldots,x_n) \in \RR^n: x_i \le M \mbox{ for all } i \in [n]\}$. Here $M$ is a parameter whose value determines the winner in the game. The problem in the test is equivalent to showing that Chooser wins if $M>2^{n-1}$, with users in the AoPS forum suspecting Chooser may still win if $M$ is polynomial in $n$~\cite{artof}.

\smallskip
This game is unusual as $K_M$ is not bounded. We remark immediately that it makes no sense to have both $v$ and $-v$ in $V$ and we redefine $V$ as the set of all vectors $v=(v_1,\ldots,v_n)$ with $v_1=1$ and $v_i\in \{1,-1\}$ for $i>1$. In fact any other choice between $v$ and $-v$ for all $\pm 1$ vectors would do.

\medskip
In Section~\ref{sec:bgame} we determine the exact smallest possible value of $M$ for which Chooser wins, which is exponential in $n$. 

\begin{thm}\label{th:bgame}
    Given the above setup, $n>1$, consider the game $G(V,K_M)$.
    \begin{enumerate}
        \item If $n$ is odd, then Chooser has a winning strategy if and only if
        \[M\geq 2^{n-2}-\frac12\binom{n-1}{\frac{n-1}{2}}.\]
        \item If $n$ is even and not a power of 2, then Chooser has a winning strategy if and only if 
        \[M\geq 2^{n-2}-\frac12\binom{n-1}{\frac{n}{2}}.\]
        \item If $n$ is a power of 2, then Chooser has a winning strategy if and only if
        \[M\geq 2^{n-2}-\frac12\binom{n-1}{\frac{n}{2}}+\frac12.\]
    \end{enumerate}
    In each case the critical value of $M$ is $2^{n-2}(1-O(n^{-1/2}))$.
\end{thm}

This theorem or rather its proof has a combinatorial application. 
\begin{thm}\label{th:comb}
If $m \ge 2$ and is not a power of 2, then the $m$-sets of a $2m$-set can be coloured Red and Blue so that 
\begin{enumerate}[label=(\roman*)]
        \item complementary $m$-sets have distinct colours,
        \item every point of the $2m$-set is contained in the same number of Red and Blue sets. 
\end{enumerate}        
        When $m$ is a power of 2, then there is again a Red-Blue colouring of the $m$-sets of a $2m$-set satisfying (i) such that for $m/2$ of the points, there are three more Red sets containing it than Blue sets, and for the remaining $3m/2$ of the points, there is one more Blue set containing it than Red sets. 
\end{thm}

We remark that when $m$ is a power of 2, this result is best possible in the sense that the number of Red and Blue sets containing any point is as equal as possible. Indeed, if $R(i),B(i)$ are the number of Red and Blue sets containing the point $i$ for $i=1,\ldots,2m$, then $\sum_i R(i)=\sum_i B(i)$ and $R(i)+B(i)=\binom{2m-1}{m}$ is odd. Furthermore, by considering what happens if we were to change colours for a complementary pair of $m$-sets, we see that $R(i)-B(i)\pmod{4}$ is independent of $i$. Thus, the smallest absolute values $R(i)-B(i)$ can take is when $R(i)-B(i)=3$ for $m/2$ different values of $i$, and $R(i)-B(i)=-1$ for the remaining $3m/2$ values of $i$.

\section{V-closed sets}\label{sec:vclos}

Although $P(V)$ seems to be the simplest $V$-closed set, there are several others. For instance, Figure 1 left shows an unbounded $V$-closed set in the plane with $V=\{e_1,e_2\}$. This $V$-closed set consists of the lattice points in the strip $\{(x_2,x_1)\in \RR^2: 0\le x_2-x_1 \le 1\}$. More generally let $S$ be a subspace of $\RR^n$, $S^{\perp}$ its complementary subspace and let $\pi$ denote the orthogonal projection $\RR^n \to S$. If $T'$ is a $\pi(V)$-closed set in $S$, then $T'+S^{\perp}$ is $V$-closed in $\RR^n$. Note however that even if condition (\ref{eq:nonpar}) holds for $V$, it may fail for $\pi(V)$. We remark further that in all of these examples the convex hull of the $V$-closed set contains a line.

\medskip
Given a convex set $C \subset \RR^n$, a point $x$ on its boundary is called \emph{exposed} if there is a supporting hyperplane $H$ to $C$ such that $\{x\}=H\cap C$. Every exposed point is an extreme point but an extreme point of $C$ is not necessarily exposed. 
We need a few facts from convex geometry. The first one is a theorem of Straszewicz~\cite{Strasz}, see also Theorem 1.4.7 in \cite{Schn}.

\begin{lem}\label{l:expos} Assume $C \subset \RR^n$ is a closed convex set and $x$ is an extreme point of $C$. Then arbitrarily close to $x$ there are exposed points of $C$.    
\end{lem}

\begin{lem}\label{l:xinT} Assume $T$ is $V$-closed and $x$ is an exposed point of $C=\cl \conv T$. Then $x \in \cl T$.
\end{lem}
{\bf Proof.}
Let $H$ be the supporting hyperplane to $C$ with outer unit normal $a$ satisfying $\{x\}=C\cap H$. So $C$ is contained in the halfspace $H^-=\{y\in \RR^n: a\cdot y\le a\cdot x\}$. The cap $C_s$ of $C$ is defined as $\{y\in C: a\cdot y\ge a\cdot x-s\}$ where $s>0$ is small. Let $r_s$ be the radius of the smallest ball with center $x$ and containing $C_s$. 

\medskip
We {\bf claim} that $C_s$ shrinks to the point $x$ as $s\to 0$. Assume on the contrary that there is a sequence $s_m \to 0$ as $m$ goes to infinity, and there is a point $x_m\in C_{s_m}$ such that $\|x_m-x\|\ge \delta$ for all $m$ with a suitable $\delta>0$. We can assume that $\|x_m-x\|=\delta$ because if not we replace $x_m$ by the point on the segment $[x,x_m]$ at distance $\delta$ from $x$. By compactness the sequence $x_m$ has a convergent subsequence with limit $x^*$, say. This point $x^*$ is in $C$, is at distance $\delta$ from $x$, and since each $x_m \in C_{s_m}$, $x^*$ is on $H$. So $x^*$ is on the boundary of $C$ and $x\ne x^*$. Consequently the segment $[x,x^*]$ lies in $C$. This is impossible as $x$ is an exposed point of $C^*$, finishing the proof of the claim.

\medskip
As $x \in C$ and $C$ is the closure of $\conv T$, the cap $C_s$ must contain a point, say $x_s$, from $T$ for every $s>0$. Then $\lim x_s=x$ and $x$ lies indeed in $\cl T$.\qed

\medskip
For a (finite) set $W\subset \RR^n$ the zonotope $Z(W)$ is the polytope defined as  
\[Z(W)=\left\{\sum_{w\in W}\alpha(w)w: \alpha(w)\in [0,1]\; \forall w\in W\right\}.
\]
It is clear that $\conv P(V)$ is the zonotope $Z(V)$ and, just like $P(V)$, $Z(V)$ is also centrally symmetric with center $g(V)$. 

\medskip
In the next section we are going to prove Theorem~\ref{th:vclos} under the extra condition that $x$ is an exposed point of $C$. This implies the theorem for extreme points by the following argument. 

\medskip
Let $x$ be an extreme point of $C$. By Lemma~\ref{l:expos} $x$ is the limit of exposed points $x_m \in C$. According to Lemma~\ref{l:xinT}   each $x_m\in T$ (so $x \in T$), and there is a translate $t_m+P(V)$ of $P(V)$ containing $x_m$ and contained in $C$. Let $H_m$ be the supporting hyperplane to $C$ at $x_m$ with outer unit normal $a_m$. Then $a_m\cdot v\ne 0$ as otherwise either $x_m+v$ or $x_m-v$ is  in $T\subset C$ and in $H_m$ contradicting the fact that $x_m$ is an exposed point. Then there is a unique vertex of $\conv P(V)$, say $p_m$, where the outer normal to $\conv T$ is $a_m$. Consequently $x_m=t_m+p_m$. 

\medskip
It is clear that $\conv P(V)$ is a zonotope with finitely many, at most $2^{|V|}$ vertices. So there is a vertex of $\conv P(V)$, say $p$, such that $p_m=p$ for infinitely many $m$. Along this subsequence $x_m=t_m+p$. So $x_m$ tends to $x$ and then $t_m$ tends to $t=x-p$ along this subsequence and $t_m+P(V)$ tends to the translate $t+P(V)$ that contains $x$ and is contained in $C$.

\bigskip
\section{Proof of Theorem~\ref{th:vclos} for exposed points}\label{sec:proof}

We split this proof into two parts: first we deal with the planar case $n=2$, and then with the general case $n>2$. In both parts we assume that $C$ is unbounded because the bounded case was proved in \cite{Bar79}. We mention, however, that the proof given below shows that in the case when $C$ is bounded, for every extreme point $x \in C$ there is a translate of $P(V)$ with $x\in t+P(V)\subset C.$ This is 
slightly stronger than the corresponding lemma in \cite{Bar79}. 

\medskip
{\bf Proof} of Theorem~\ref{th:vclos} in the planar case ($n=2$). We begin with a simple fact that is mentioned also in \cite{Bar79}. We omit the proof.

\begin{fact}\label{fact:lines} If $L$ is a supporting line to $C$ and is parallel with some $v \in V$, then $L\cap C$ is either a segment whose length is at least the length of $v$ or is a halfline.
\end{fact}

The recession cone $C^*$ of $C$ is defined as $C^*=\{z\in \RR^n:z+C\subset C\}$. As $C$ is unbounded, its recession cone contains a halfline. We may assume that this halfline is the vertical one $F=\{(0,y)\in \RR^2: y\ge 0\}$. Then $z \in C$ implies that $z+C^*\subset C$. We assume that $F$ is in the interior of $C^*$ and return later to the case when $F=C^*$. 

\medskip
We have to fix some notation. Let $U\subset V$ be the set of those $v\in V$ for which there is a line $L(v)$ in the direction $v$ supporting $C$. For $v \in U$, $L(v)\cap C$ is a segment $[b,c]$.  This happens for all $v \in U$ except possibly for the $v \in U$ with the smallest and largest slope and then $L(v)\cap C$ may be a halfline, to be denoted by $[b,c]$. These segments/halflines are on the boundary of $C$ and the exposed point $x\in C$ lies between two such consecutive segments.

\begin{figure}[h!]
\centering
\includegraphics[scale=0.9]{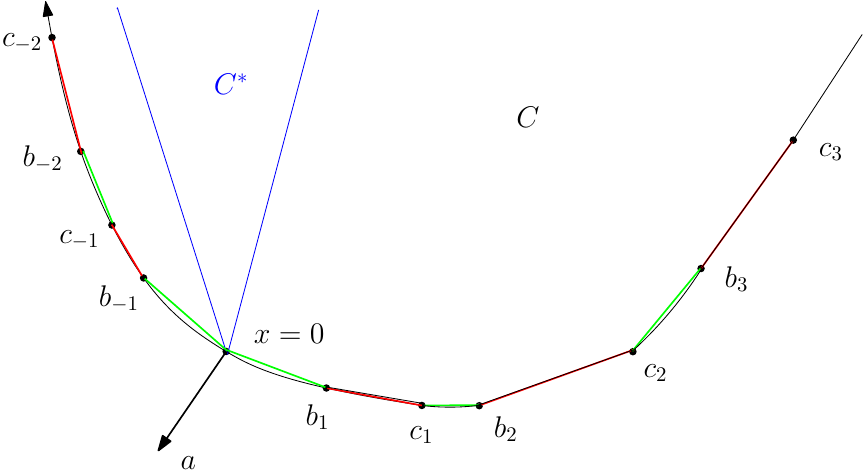}
\caption{$C$, $C^*$ and the lower boundary of $\conv D$.}
\label{fig:CandC}
\end{figure}

\medskip
Let $a$ be the outer unit normal to $C$ at $x$. We assume that $x=0$ (which can be reached by a translation), and assume further that $a\cdot v>0$ for every $v\in V$ which can be achieved by changing the sign of some vectors in $V$. For $z\in \RR^2$ we set $z=(z^x,z^y)$. Define $V^0=V\cap C^*$, 
$V^+=\{v\in V\setminus V^0: v^x>0\}$, $V^-=\{v\in V\setminus V^0: v^x<0\}$. 
The order of the vectors in $V^+$ by increasing slope is $v_1,v_2,\ldots,v_k$, and order of the ones in $V^-$ 
by decreasing slope is $v_{-1},v_{-2},\ldots,v_{-h}$. Then $L(v_i)\cap C=[b_i,c_i]$ and $b_i^x< c_i^x$ and $c_i-b_i=\beta_i v_i$ with $\beta_i\ge 1$ because of Fact~\ref{fact:lines}. Similarly $L(v_{-i})\cap C=[b_{-i},c_{-i}]$ and $b_{-i}^x >c_{-i}^x$ and $c_{-i}-b_{-i}=\beta_{-i} v_{-i}$ with $\beta_{-i}\ge 1$, see Figure~\ref{fig:CandC}. 

\medskip

Define $U^*=V^+\cup V^-$ and $D$ as the the convex hull of the point $x=0$ and of the segments (possibly halflines) $\setcond{L(v_i)\cap C}{i\in [k]}$ and $\setcond{L(v_{-j})\cap C}{j\in [h]}$. In Figure~\ref{fig:CandC} the lower boundary of $\conv D$ is coloured red and green. Set 
\[
W=\{b_1,c_1-b_1,b_2-c_1,c_2-b_2,\ldots, c_k-b_k\}\cup \{b_{-1},c_{-1}-b_{-1},b_{-2}-c_{-1},c_{-2}-b_{-2},\ldots, c_{-h}-b_{-h}\},
\]

\begin{claim}\label{cl:UW} $Z(U^*) \subset Z(W)$ and $Z(W) \subset \conv D+C^*$ and $\conv D+C^*\subset C$.
\end{claim}

{\bf Proof.} The first inclusion is very simple: $z\in Z(U^*)$ is of the form $\sum_{u\in U^*}\alpha(u)u$, and every $u \in U^*$ is $\beta(w)w$ for some $w\in W$ with $0<\beta(w)\le 1$. So $z$ is the sum of $(\alpha(u)\beta(w))w$ with the corresponding $w$ and the coefficient $(\alpha(u)\beta(w))\le 1$. Consequently this sum is in $Z(W)$.

\medskip
For the second inclusion $z\in Z(W)$ can be written as $z=\sum_{w\in W} \alpha(w)w$ where $\alpha(w)\in[0,1]$ and assume that $z^x\ge 0$. Then $z^x \in [b_i^x,c_i^x]$ or $[c_{i-1}^x,b_i^x]$ (possibly with $c_0=0$). In the first case 
$z^x=b_1^x+(c_1-b_1)^x+\ldots +(b_i-c_{i-1})^x+ \gamma(c_i-b_i)^x$ where $\gamma \in [0,1]$. The point $z^*:=b_1+(c_1-b_1)+\ldots +(b_i-c_{i-1})+ \gamma(c_i-b_i)\in Z(W)$ is on the boundary of $Z(W)$. Consequently $z\in z^*+F$. The same argument works when $z^x<0$. 

\medskip
The last inclusion follows as $D \subset C$ and $C^*$ is the recession cone of $C.$ This finishes the proof of the claim.

\medskip
Finally $V = V^+\cup V^-\cup V^0$ and $P(V)\subset Z(V) \subset Z(U^*) +C^*$ and $P(V) \subset C$ follows. 

\medskip
The case when $F=C^*$ is simple and is left to the reader. The same applies to the case when a generator of the recession cone coincides with some $v \in $V.
\qed

\bigskip
{\bf Proof} of Theorem~\ref{th:vclos} for the general case $n>2$. 

\medskip
Assume $x\in T$ is an exposed point of $C$. There is a hyperplane $H$ supporting $C$ with outer unit normal $a$ such that $\{x\}=H \cap C$. 
We check the same way as above that $a\cdot v\ne 0$ for every $v \in V$. Again $\conv P(V)$ is the zonotope generated by the vectors in $V$. Its vertex where the outer normal to $\conv P(V)$ is $a$ is uniquely determined. Let $p$ be this vertex.

\smallskip
Let $b$ be a unit vector orthogonal to $a$ and let $S$ be the two-dimensional plane containing $x$, $x+a$ and $x+b$. Denote by $\pi$ the orthogonal projection $\RR^n \to S$. We want to choose the vector $b$ so that the set $\pi(V)$ still satisfies condition (\ref{eq:nonpar}). The fact that $\pi(u)$ and $\pi(v)$ are parallel means (via a simple computation) that $b\cdot [(a\cdot v)u-(a\cdot u)v]=0$. Here $(a\cdot v)u-(a\cdot u)v\ne 0$ because  $a\cdot v,a\cdot u\ne 0$ and the vectors $u,v$ are not parallel. Moreover $(a\cdot v)u-(a\cdot u)v$ is not parallel with $a$ because if $(a\cdot v)u-(a\cdot u)v=\lambda a$, then $\lambda a^2=a\cdot[(a\cdot v)u-(a\cdot u)v]=0$ implies $\lambda=0$. Thus $b\cdot [(a\cdot v)u-(a\cdot u)v]=0$ is a linear condition on $b$ and there are finitely many such conditions. Then we can choose $b$ so that $\pi(V)$ satisfies condition (\ref{eq:nonpar}). In fact, choosing $b$ randomly (with $a\cdot b=0$) suffices with probability one.

\smallskip
Our next target is to show that there is a translate of $\pi(P(V))=P(\pi(V))$ that contains $\pi(x)=x$ and is contained in $\pi(C)$. This is exactly the two-dimensional case in the plane $S$ that we just proved. Let this translate be $t+\pi(P(V))$. Moreover, the translate $t+\pi(P(V))$ is the projection by $\pi$ of $x-p+P(V)$ because the outer normal to  $P(V)$ at vertex $p$ and the outer normal to $C$ at $x$ is the same vector $a$. This means that for a randomly chosen vector $b$ (with $a\cdot b=0$) the projection onto $S$ of $x-p+P(V)$ contains $x$ and is contained in $\pi(C)$. 

\begin{claim} The translate $x-p+P(V)$ contains $x$ and is contained in $C$. 
\end{claim}

{\bf Proof.} We only have to check that $x-p+P(V)\subset C$. Assume that this is not the case. Then $x-p+q \notin C$ for some vertex $q$ of $\conv P(V)$. Then there is a separating hyperplane that contains $x-p+q$ and is disjoint from $C$. As the separation is strict, the  normal, $a^*$ say, of this hyperplane can be chosen with some freedom. 

\medskip
Let $b$ be the vector orthogonal to $a$ in the 2-dimensional plane $S$ containing $x, x+a$ and $x+a^*$. Again, let $\pi$ denote the orthogonal projection $\RR^n \to S$. We can choose $a^*$ in such a way that the vectors in $\pi(V)$ satisfy condition (\ref{eq:nonpar}). The planar case shows that the projection of $x-p+P(V)$ is contained in $\pi(C)$. In particular $\pi(x-p+q) \in \pi(C)$. So the hyperplane with normal $a^*$ can not separate $x-p+q$ and $C$.

\medskip
This finishes the proof of the claim and of Theorem~\ref{th:vclos}\qed

\section{A special balancing game}\label{sec:bgame}

In this section, we analyze the game $G(V,K_M)$ for $V=\setcond{v\in \{-1,1\}^n}{v_1=1}$, where $K_M$ is defined as $\{(x_1,\ldots,x_n) \in \RR^n: x_i \le M \mbox{ for all } i \in [n]\}$, and prove Theorem~\ref{th:bgame}. First, we note the following fact which distinguishes the case when $n$ is a power of 2.
\begin{lem}\label{lem:V0size}
    Let $n$ be even. Then $\frac12 \binom{n}{n/2}$ is an odd integer if $n$ is a power of 2, and an even integer otherwise.
\end{lem}
{\bf Proof.}
    Let $k$ be the largest integer such that $n\geq 2^k$. For an integer $x$, denote by $\nu_2(x)$ the largest integer $d$ such that $2^d\mid x$. Then $\nu_2(n!)=\floor{n/2}+\floor{n/4}+\cdots+\floor{n/2^k}$. Thus,
    \begin{align*}
    &\nu_2\paren{\binom{n}{n/2}}=\nu_2(n)-2\nu_2(n/2)\\
    &=(\floor{n/2}-2\floor{n/4})+(\floor{n/4}-2\floor{n/8})+\cdots+(\floor{n/2^{k-1}}-2\floor{n/2^k})+\floor{n/2^k}\geq 1,
    \end{align*}
    with equality if and only if $\floor{n/2^i}-2\floor{n/2^{i+1}}=0$ for $i=1,\ldots,k-1$. This is equivalent to $n$ being a power of 2.
\qed

\medskip
We begin with the strategy for Pusher.

\subsection{Pusher's strategy}

By Theorem~\ref{th:determined}, we only need to show that if there exists a $V$-closed set $T$ with $0\in T\subseteq K_M$, then $M$ must be at least the values given in Theorem~\ref{th:bgame}.

\bigskip
Suppose there exists a $V$-closed set $T$ with $0\in T\subseteq K_M$, for some $M\in \ZZ$. The intersection of $T$ with $\ZZ^n$ is also $V$-closed, so we may assume that $T\subset \ZZ^n$, so that $T$ is closed (even discrete). Let $C=\cl\conv T$. Since $C\subseteq K_M$, it contains no line. Let $\one=(1,\ldots,1)\in \ZZ^n$ be the all-ones vector. Let $H$ be the supporting hyperplane of $C$ perpendicular to $\one$, so $H$ is of the form $\setcond{x\in \RR^n}{x\cdot \one=s}$ where $s\geq 0$ is the smallest real number such that $C$ lies in the half-space $\setcond{x\in \RR^n}{x\cdot \one\leq s}$. Here $s\ge 0$ follows from $0\in T\subset C$.

\medskip
Let $C'=C\cap H$. Since $C$ is closed and convex and $H$ is a supporting hyperplane, $C'$ is non-empty and convex. Since $C'\subseteq H\cap K_M$ is bounded, there is an extreme point $x$ of $C'$, hence is also an extreme point of $C$. By Theorem~\ref{th:vclos}, there is a translation vector $t\in \RR^n$ such that $x\in t+P(V)\subset C\subset K_M$. Since $x$ is an extreme point, we have $x\in T\subset \ZZ^n$. 

\medskip 
We now have a translate $t+P(V)$ which intersects the hyperplane $H$ and lies completely in $K_M$. We note that the smallest cube containing $P(V)$ is a translate of $[0,2^{n-1}]^n$. This follows from the fact (which is easy to check) that the width of $P(V)$ in direction $e_i$ is $2^{n-1}$ for every $i \in [n]$. Then $P(V)$ touches all $2n$ facets of this cube. It is then clear that if a translate of $P(V)$ is contained in $K_M$, then the smallest cube containing this translate is also contained in $K_M$. Thus the furthest a translate of $P(V)$ can be along the $\one$ direction is when it is touching all $n$ facets of $K_M$. Thus, computing the translate of $P(V)$ touching all $n$ faces of $K_M$ is a simple task. However, we will show that this translate is not far enough along $\one$ to touch $H$ for our choice of $M$. 

\medskip
Since $t+P(V)\in K_M$, for $i=1,\ldots,n$, we have $t_i+\max_{u\in P(V)} u_i\leq M$. Recall that $V=\setcond{v\in \{-1,1\}^n}{v_1=1}$. For $i=1$, this implies that $t_1\leq M-|V|=M-2^{n-1}$. For $i>1$, we have $t_i\leq M-2^{n-2}$. We now compute
\[r:=\max_{u\in P(V)} (u\cdot \one)=\sum_{v\in V,v\cdot \one>0} v\cdot \one.\]
When $n$ is odd, 
\[r=\sum_{j=0}^{\frac{n-1}{2}} (n-2j)\binom{n-1}{j}=\frac{n}{2}\binom{n-1}{\frac{n-1}{2}}+2^{n-2}.\]
When $n$ is even,
\[r=\sum_{j=0}^{\frac{n}{2}-1} (n-2j)\binom{n-1}{j}=\frac{n}{4}\binom{n}{n/2}+2^{n-2}=\frac{n}{2}\binom{n-1}{n/2}+2^{n-2}.\]
Since $t+P(V)$ intersects $H$, there exists $u\in P(V)$ such that $(t+u)\cdot \one=s\geq 0$. On the other hand,
\[(t+u)\cdot \one\leq nM-(n+1)2^{n-2}+u\cdot \one\leq nM-(n+1)2^{n-2}+r.\]
Thus, we have
\[M\geq 2^{n-2}+\frac{1}{n}(2^{n-2}-r)=\begin{cases}
    2^{n-2}-\frac12 \binom{n-1}{\frac{n-1}{2}} & \text{if $n$ is odd,}\\
    2^{n-2}-\frac{1}{2}\binom{n-1}{n/2} & \text{if $n$ is even.}
\end{cases}\]
Since $M\in \ZZ$, we get the desired bounds in Theorem~\ref{th:bgame} for all $n$, noting that when $n$ is a power of 2, by Lemma~\ref{lem:V0size}, $\binom{n-1}{n/2}=\frac12\binom{n}{n/2}$ is odd. Then, $M\geq \ceil{2^{n-2}-\frac{1}{2}\binom{n-1}{n/2}}=2^{n-2}-\frac{1}{2}\binom{n-1}{n/2}+\frac12$.

\subsection{Chooser's strategy}

Next, we construct Chooser's winning strategy in Theorem~\ref{th:bgame}. From Theorem~\ref{th:determined}, her strategy is to simply construct a $V$-closed set $T$ such that $0\in T\subset K_M$. 
Note that it suffices to find a translation vector $t\in \RR^n$ such that $0\in t+P(V)\subset K_M$. Our translate can be described as follows. First, center $P(V)$ at the origin, then translate it along the direction $-\one$ for as long as the convex hull of the translated $P(V)$ contains the origin. If $n$ is not a power of 2, then this translate of $P(V)$ will contain the origin. If $n$ is a power of 2, then this translate of $P(V)$ will not contain the origin for parity reasons, so we need a minor adjustment.

\medskip
Write $\Sigma(V)=\sum_{v\in V} v=2^{n-1}e_1$. Recall that the center of $P(V)$ is $g(V)=\frac12 \Sigma (V)$. We consider the translate  

\begin{equation}\label{eq:PVcenter}
    P(V)-g(V)=\setcond{\frac12 \sum_{v\in V} \varepsilon_v v}{\varepsilon_v\in \{-1,1\}}. 
\end{equation}
The center of $P(V)-g(V)$ is the origin. We now deal with the even and odd cases in Theorem~\ref{th:bgame} separately. 

\medskip
{\bf Case 1.} When $n$ is odd. Set $M=2^{n-2}-\frac12 \binom{n-1}{\frac{n-1}{2}}$. We claim that 
\[0\in P(V)-g(V)-\frac12\binom{n-1}{\frac{n-1}{2}}\one\subset K_M.\]
These two inclusions follows from (\ref{eq:PVcenter}) and the next two claims.
\begin{claim}
    There exists $\varepsilon_v\in \{-1,1\}$ for each $v\in V$ such that $\sum_{v\in V} \varepsilon_v v=\binom{n-1}{\frac{n-1}{2}}\one$.
\end{claim}
{\bf Proof.}
    For each $v\in V$, set
    \[\varepsilon_v=\begin{cases}
        1 & \text{if } \sum_i v_i>0,\\
        -1 & \text{if } \sum_i v_i<0,
    \end{cases}\]
    noting that the sum cannot be 0 since $n$ is odd. For each index $i$, the $i$-th coordinate of $\sum_{v\in V} \varepsilon_v v$ is $A-B$, where
    \begin{align*}
        A &= \#\setcond{v\in V}{v_i=1,\ \textstyle\sum_j v_j>0}=\sum_{k=\frac{n-1}{2}}^{n-1}\binom{n-1}{k}\\
        B &= \#\setcond{v\in V}{v_i=-1,\ \textstyle\sum_j v_j>0}=\sum_{k=\frac{n-1}{2}+1}^{n-1}\binom{n-1}{k}.
    \end{align*}
    Thus, $A-B=\binom{n-1}{\frac{n-1}{2}}$ and hence $\sum_{v\in V} \varepsilon_v v=\binom{n-1}{\frac{n-1}{2}}\one$.
\qed

\begin{claim}\label{claim:subsetKM}
    For any $\varepsilon_v\in \{-1,1\}$ for each $v\in V$, we have $\frac12 \sum_{v\in V} \varepsilon_v v-\frac12 \binom{n-1}{\frac{n-1}{2}}\one\subset K_M$.
\end{claim}
{\bf Proof.}
    It suffices to show that for each index $i$, the $i$-coordinate of $\frac12 \sum_{v\in V} \varepsilon_v v$ is at most $2^{n-2}$. Indeed, this is true since each $\varepsilon_v v$ contributes at most 1 to coordinate $i$ and there are exactly $2^{n-1}$ vectors in $V$.
\qed

\bigskip
{\bf Case 2.} When $n$ is even. This case is more complex due to the existence of the ``middle layer'', $V_0=\setcond{v\in V}{\sum_i v_i=0}$. Let $L\subset \ZZ^n$ be the $(n-1)$-dimensional lattice generated by $V_0$. Then $P(V_0)$ lies in $L$, which can also be expressed as 
\[L=\setcond{u\in \ZZ^n}{\sum_i u_i=0,\quad  \{u_1,\ldots,u_n\} \text{ consists of integers that are all even or all odd}}.
\]
Alternatively, $L$ is generated by the vectors $\{2(e_1-e_2),2(e_1-e_3),\ldots,2(e_1-e_n),w\}$, for any $w\in V_0$. 

\medskip
We show that $P(V_0)$ contains every vector in $L$ near its center $g(V_0)$.

\begin{lem}\label{lem:PV0ball}
    Let $R$ be a positive integer such that 
    \begin{equation}\label{eqn:Rbound}
        \binom{n-2}{\frac{n-2}{2}}>4R(n-1).
    \end{equation}
    Then for any $u\in L$ with $\norm{u-g(V_0)}_{\infty} \le R-\frac{n}{2}$, we have $u\in P(V_0)$. 
\end{lem}
{\bf Proof.}
    Recall that $P(V_0)-g(V_0)$ consists of vectors of the form $\frac12\sum_{v\in V_0} \varepsilon_v v$. For $r>0$, write $B_r(v)=\setcond{x\in \RR^n}{\norm{x-v}_2\leq r}$ for the ball of radius $r$ centered at $v\in \RR^n$. The proof strategy is to find $U\subseteq V_0$ such that $P(U)$ attains every vector in $L$ up to radius $R$ from its center. Then, using the remaining vectors $V_0\setminus U$, we find an element $x\in P(V_0\setminus U)$ that is close to its center. Combining, we can attain from $P(V_0)$ every vector in $L$ close to its center.

\smallskip
    It will be convenient to define $V_0^{\pm}=V_0\cup -V_0=\setcond{v\in \pm V}{\sum_i v_i=0}$ and the equivalence relation $v\sim v'$ if $v=\pm v'$. We say that $v,v'$ are \emph{$(\sim)$-distinct} if $v\nsim v'$.
 
\smallskip   
    We shall find $2R(n-1)$ $(\sim)$-distinct vectors $v_i^{(j)+},v_i^{(j)-}\in V_0^{\pm}$, for $i=2,\ldots,n$ and $j=1,\ldots,R$, such that $v_i^{(j)+}-v_i^{(j)-}=2(e_1-e_i)$. We do this greedily, by finding a pair $v_i^{(j)+},v_i^{(j)-}$ in any ordering of the pairs of indices $(i,j)$, that is distinct from all the other vectors we have found so far. 

\smallskip
    Suppose we have already found $m$ pairs, and it is time to find $v_i^{(j)+},v_i^{(j)-}$ for the $(m+1)$-th pair of indices $(i,j)$. Choose some $v_i^{(j)+}$ so that $(v_i^{(j)+})_1=1$ and $(v_i^{(j)+})_i=-1$, then set $v_i^{(j)-}$ to be identical to $v_i^{(j)+}$, except its first and $i$-th coordinate have flipped signs. This satisfies $v_i^{(j)+}-v_i^{(j)-}=2(e_1-e_i)$. Note that there are $\binom{n-2}{\frac{n-2}{2}}$ ways to pick $v_i^{(j)+}$, at most $4m$ of which will result in a collision between $v_i^{(j)+}$ or $v_i^{(j)-}$ with one of the $2m$ vectors we previously found. By the assumption on $R$, 
    \[\binom{n-2}{\frac{n-2}{2}}>4R(n-1)\geq 4m,\]
    thus this greedy process will be successful in finding all $2R(n-1)$ vectors.

\smallskip
    We now have $2R(n-1)$ $(\sim)$-distinct vectors $v_i^{(j)+},v_i^{(j)-}\in V_0^{\pm}$ satisfying $v_i^{(j)+}-v_i^{(j)-}=2(e_1-e_i)$ for $i=2,\ldots,n$ and $j=1,\ldots,R$. Pick any $w\in V_0^{\pm}$ distinct from these, then set $U:=\set{v_2^{(1)+},v_2^{(1)-},\ldots,v_n^{(R)+},v_n^{(R)-},w}$. 

\smallskip
    The set $P(U)-g(U)$ consists of vectors of the form $\frac12 \sum_{u\in U}\varepsilon_u u$. To exhibit many vectors from this set, we assert that we will use each pair $v_i^{(j)+},v_i^{(j)-}$ with opposite signs, that is, either $v_i^{(j)+}-v_i^{(j)-}=2(e_1-e_i)$ or $v_i^{(j)-}-v_i^{(j)+}=2(e_i-e_1)$. Then, $P(U)-g(U)$ contains all vectors of the form
    \begin{equation}\label{eqn:PUform}
        a_2(e_1-e_2)+a_3(e_1-e_3)+\cdots+a_n(e_1-e_n)\pm\frac12 w,
    \end{equation}
    for any integer $-R\leq a_i\leq R$ with $a_i\equiv R\pmod{2}$. Using the next two claims, we show that all vectors in $(L-g(U))\cap B_R(0)$ are expressible in this form. 

    \begin{claim} \label{claim:L-12SigmaU}
        $L-g(U)$ consists of all vectors of the form
        \[b_2(e_1-e_2)+b_3(e_1-e_3)+\cdots+b_n(e_1-e_n)\pm\frac12 w,\]
        where $b_i\in \ZZ$, $b_i\equiv R\pmod{2}$. 
    \end{claim}
    {\bf Proof.}
        We know that $L$ consists of all vectors of the form
        \[b_2(e_1-e_2)+b_3(e_1-e_3)+\cdots+b_n(e_1-e_n)+\varepsilon w,\quad b_i\equiv 0\pmod{2}, \quad \varepsilon\in \{0,1\}.\]
        Set $z:=R(e_1-e_2)+R(e_1-e_3)+\cdots+R(e_1-e_n)-\frac12 w$. Since $v_i^{(j)+}-v_i^{(j)-}=2(e_1-e_i)$ and $V_0\subset L$, we have
        \[z=\sum_{i=2}^n \sum_{j=1}^R \frac12(v_i^{(j)+}-v_i^{(j)-})-\frac12 w\in L-g(U)\]
        Thus, $z+L=L-g(U)$ consists of all vectors of the form
        \[b_2(e_1-e_2)+b_3(e_1-e_3)+\cdots+b_n(e_1-e_n)\pm\frac12 w,\quad b_i\equiv R\pmod{2}\]
    since $V_0\subset L$.
    \qed

    \begin{claim} \label{claim:LcapBsubsetPU}
        For any $x\in L$ with $\norm{x-g(U)}_\infty \leq R$, we have $x\in P(U)$.
    \end{claim}
    {\bf Proof.}
        From Claim~\ref{claim:L-12SigmaU}, we may write $x-g(U)=a_2(e_1-e_2)+a_3(e_1-e_3)+\cdots+a_n(e_1-e_n)\pm\frac12 w$. Since $\norm{x-g(U)}_\infty\leq R$, by looking at the $i$-th coordinate of $x$ for $i=2,\ldots,n$, we have $|a_i|\leq R+\frac12$, so in fact we have $|a_i|\leq R$. Therefore, $x-g(U)$ is of the form (\ref{eqn:PUform}). This implies that $x\in P(U)$.
    \qed

\smallskip
Each $v\in V_0\setminus U$ has norm $\norm{v}_2=\sqrt{n}$. Theorem 2.2 in \cite{beck} (or Lemma 4.1 in \cite{Bar08}) shows that there exists for each $v\in V_0\setminus U$, a sign $\varepsilon_v\in \set{-1,1}$ such that for $x:=\sum_{v\in V_0\setminus U} \varepsilon_v v$, we have $\norm{x}_2\leq n$ and $\frac12(x+\Sigma(V_0\setminus U))=\frac 12 x + g(V_0\setminus U) \in P(V_0\setminus U)\subset L$. Then of course $\norm{x}_{\infty}\le n$.

\medskip
    Given any $u\in L$ with $\norm{u-g(V_0)}_{\infty} \le R-\frac{n}{2}$, we have $\norm{u-(\frac12 x+g(V_0)}_{\infty}\le R$ and thus $u-(\frac12 x+g(V_0\setminus U))\in L\cap B_R(g(U))$. By Claim~\ref{claim:LcapBsubsetPU}, $u-(\frac12 x+g(V_0\setminus U))\in P(U)$. Since $\frac12 x+g(V_0(W)\setminus U))\in P(V_0\setminus U)$, we have $u\in P(V_0)$ as required.
\qed

\bigskip

\smallskip
In the next lemma, we show that $P(V_0)$ contains its center $g(V_0)$ if $n$ is not a power of 2. In the case when $n$ is a power of 2, $P(V_0)$ contains the closest possible point to its center $g(V_0)$.
\begin{lem} \label{lem:nevencenter}
    Suppose $n$ is even. 
    \begin{enumerate}
        \item If $n$ is not a power of 2, then $g(V_0)\in P(V_0)$.
        \item If $n$ is a power of 2 and $n>2$, let $w\in \ZZ^n$ be given by
        \[w_i=\begin{cases}
            1 & \text{if } i\not\equiv 0\pmod{4},\\
            -3 & \text{if } i\equiv 0\pmod{4}.
        \end{cases}\]
        Then, $g(V_0)-\frac12 w\in P(V_0)$.
    \end{enumerate}
\end{lem}
{\bf Proof.}
    Recall that $L$ is the lattice generated by $V_0$. 
    \begin{enumerate}
        \item We first show that $g(V_0)\in L$. Indeed, its first coordinate is $(g(V_0))_1=\frac12 \Sigma (V_0)_1=\frac12 |V_0|$, an integer by Lemma~\ref{lem:V0size}. In fact $|V_0|={n \choose \frac n2}$ but we don't need this explicitely. By symmetry and the fact that the sum of the coordinates of $g(V_0)$ is 0, we have for $i=2,\ldots,n$, $(g(V_0))_i=-\frac{1}{2(n-1)} |V_0|$, also an integer. Since $n$ is even, all the coordinates have the same parity, thus $g(V_0)\in L$.
 
        \smallskip   
        By Lemma~\ref{lem:PV0ball} with $R=n/2$ and $u=g(V_0)$, we have $g(V_0)\in P(V_0)$, provided that (\ref{eqn:Rbound}) holds. It can be shown that (\ref{eqn:Rbound}) holds for $n\geq 14$. For even $n<14$, we explicitly construct the signs $\varepsilon_v$ in the Appendix, see Lemma~\ref{lem:nevensmall}.

        \item Let $u=g(V_0)-\frac12 w$. As above, we first show that $u\in L$. Indeed, it is straightforward to check that $u_1=\frac12 |V_0|-\frac12 w_1$ and $u_i=-\frac{1}{2(n-1)}|V_0|-\frac12 w_i$ for $i\geq 2$. Furthermore, $\sum_i u_i=0$ and $u_1\equiv u_2\equiv \cdots\equiv u_n\pmod{2}$, thus $u\in L$. 
        
        \smallskip
        Applying Lemma~\ref{lem:PV0ball} with $R=n/2+3/2$ and $u\in L$, checking that (\ref{eqn:Rbound}) holds for $n\geq 16$, we have that $u\in P(V_0)$, as required. We leave the remaining small cases $n<16$ in the Appendix, see Lemma~\ref{lem:npow2small}.
    \end{enumerate}
\qed

\smallskip
The next two lemmas give Chooser's strategy for when $n$ is not a power of 2, and when $n$ is a power of 2. Together with the $n$ odd case, this completes the Chooser's strategy portion of Theorem~\ref{th:bgame}.
\begin{lem}\label{lem:neven}
    Suppose $n$ is even and not a power of 2. Let $M= 2^{n-2}-\frac12\binom{n-1}{n/2}$. Then
    \[0\in P(V)-g(V)-\frac12\binom{n-1}{n/2}\one\subset K_M.\]
\end{lem}
{\bf Proof.}
    The proof of the second inclusion is identical to the proof of Claim~\ref{claim:subsetKM}.

\smallskip
    For the first inclusion, from (\ref{eq:PVcenter}), it suffices to find $\varepsilon_v\in \{-1,1\}$ for each $v\in V$ such that $\frac12\sum_{v\in V} \varepsilon_v v=\frac14\binom{n}{n/2}\one$. We first set for $v\in V\setminus V_0$,
    \[\varepsilon_v=\begin{cases}
        1 & \text{if }\sum_i v_i>0,\\
        -1 & \text{if }\sum_i v_i<0.
    \end{cases}\]
    Then it is not hard to see that $\sum_{v\in V\setminus V_0} \varepsilon_v v=\frac12\binom{n}{n/2}\one=
    \binom{n-1}{n/2}\one$. Thus, it suffices to find $\varepsilon_v\in \{-1,1\}$ for each $v\in V_0$ such that $\frac12\sum_{v\in V_0} \varepsilon_v v=0$. Equivalently, this says that $g(V_0)=\frac12 \Sigma (V_0)\in P(V_0)$, which is Lemma~\ref{lem:nevencenter}(1).
\qed

\begin{lem} \label{lem:npow2}
    Suppose $n=2^k$ for some $k\geq 1$. Let $M= 2^{n-2}-\frac12\binom{n-1}{n/2}+\frac12$. Then there exists $w\in \ZZ^n$ with $w_i\leq 1$ for all $i$, such that
    \[0\in P(V)-g(V)-\frac12\binom{n-1}{n/2}\one+\frac12 w\subset K_M.\]
\end{lem}
{\bf Proof.}
    The case $n=2$ can be directly verified, with $w_1=1,w_2=-1$. Assume that $k\geq 2$, so that $4\mid n$. Let $w\in \ZZ^n$ be given by
    \[w_i=\begin{cases}
        1 & \text{if } i\not\equiv 0\pmod{4},\\
        -3 & \text{if } i\equiv 0\pmod{4}.
    \end{cases}\]
    The proof of the second inclusion is identical to the proof of Claim~\ref{claim:subsetKM}. For the first inclusion, following the same arguments as in Lemma~\ref{lem:neven}, it suffices to show that $u:=g(V_0)-\frac12 w\in P(V_0)$, which is just Lemma~\ref{lem:nevencenter}(2).
\qed

\bigskip
{\bf Remark.} Our methods also extend to determining the exact threshold $M$ in other variants of the balancing game $G(V,K_M)$. For instance, one can consider the case where $V=\{-1,0,1\}^n$, or where $K_M=\setcond{(x_1,\ldots,x_n)\in \RR^n}{|x_i|\leq M \text{ for all }i\in [n]}$, or both. We leave the determination of the precise threshold for these variants to the interested reader.

\section{A combinatorial consequence, proof of Theorem~\ref{th:comb}}

Let the $2m$ set be $[2m]$. With an $m$-set $A$ we associate the vector $I(A) \in \RR^{2m}$ whose $i$th component is $1$ if $i\in A$ and $-1$ if $i \notin A$. This creates a set of $N:={2m \choose m}$ vectors. Note that $I(A)=-I(A^c)$ where $A^c$ is the complement $m$-set of $A$. Keep just one of these complementary pairs and make an $2m \times N/2$ matrix $M$ with these vectors as columns. Suppose for concreteness that we keep every $I(A)$ with $1 \in A$. 

\medskip
Assume now that there is a $1,-1$ vector $\eps \in \RR^{N/2}$ with $M\eps=0$, and let $\eps_A$ be the coordinate corresponding to column $I(A)$ of $M$.  Colour $A$ Red if $\eps_A=1$ and colour it Blue if $\eps_A=-1$. Then condition (i) of the theorem is satisfied, and so is condition (ii) because $\sum \eps_AI(A)=0$.

\medskip
The existence of such an $\eps$ is guaranteed by Lemma~\ref{lem:nevencenter}(1) in the proof of Theorem~\ref{th:bgame}.  

\medskip
Such an $\eps\in \RR^{N/2}$ does not exist when $m$ is a power of $2$ because $N/2$ is odd. However, Lemma~\ref{lem:nevencenter}(2) shows there is an $\eps\in \{-1,1\}^{N/2}$ such that every coordinate of $M\eps$ is $+3$ or $-1$. 
\qed

\bigskip
{\bf Remark.} Our construction for Theorem~\ref{th:comb} is based on a greedy algorithm and therefore is not explicit. It is conceivable that a direct coloring exists, but we have not been able to find one.

\clearpage
\appendix

\section{Appendix}

In this section, we verify Lemmas~\ref{lem:neven} and \ref{lem:npow2} for small $n$. In the proofs of the lemmas, they reduce to finding signs $\varepsilon_v$ for each $v\in V_0$. Recall that $V_0=\setcond{v\in \{-1,1\}^n}{v_1=1, \sum_i v_i=0}$. 
\begin{lem} \label{lem:nevensmall}
    For $n\in \{6,10,12\}$, there exist signs $\varepsilon_v\in \{-1,1\}$ for each $v\in V_0$ such that $\sum_{v\in V_0} \varepsilon_v v=0$.
\end{lem}
{\bf Proof.}
    It will be more convenient to work with $V_0^{\pm}=V_0\cup -V_0$. The lemma is then equivalent to finding $\varepsilon_v\in \{-1,1\}$ for each $v\in V_0^{\pm}$ such that $\varepsilon_{-v}=-\varepsilon_v$ and $\sum_{v\in V_0^{\pm}} \varepsilon_v v=0$. This can be encoded as a SAT problem, however the size of $V_0^{\pm}$ is too large to be practical for a naive encoding. Instead, we will assign signs $\varepsilon_v$ to most of the vectors in $V_0^{\pm}$ which sum to 0, then use a SAT solver to pick signs for the remaining vectors.
    
    Let $\rho:\ZZ^n\to \ZZ^n$ be the coordinate rotation $\rho(v)=(v_2,v_3,\ldots,v_n,v_1)$. Then $\rho$ acts on $V_0^{\pm}$ and decomposes it into orbits. Note that the sum of all the vectors in each orbit is 0. Thus for each $v\in V_0^{\pm}$, if $v,-v$ are in different orbits, then we assign $\varepsilon_u=1$ for each $u$ in the orbit of $v$, and $\varepsilon_{-u}=-1$. The sum of $\varepsilon_u u$ among all such orbits is 0, so it remains to pick signs for $v\in V_0^{\pm}$ where $v,-v$ belong in the same orbit.

    We enumerate all such vectors $v$ and use a SAT solver to pick a satisfiable set of signs $\varepsilon_v$. The list of all such vectors $v$ and their corresponding signs $\varepsilon_v$ is listed below for each of the three cases $n=6,10,12$. Here, the vectors are written as a binary string, with 0 in place of -1 for brevity.
    
    \medskip
    {\bf Case $n=6$.}
    
    \medskip
    \makebox[\textwidth][c]{
    \begin{minipage}[t]{0.2\textwidth}
    \centering
    \begin{tabular}{c|c}
        $\varepsilon_v$ & $v$ \\ \hline
        + & 000111\\
        - & 001110\\
        + & 011100\\
        - & 100011\\
    \end{tabular}
    \end{minipage}
    \begin{minipage}[t]{0.2\textwidth}
    \centering
    \begin{tabular}{c|c}
        $\varepsilon_v$ & $v$ \\ \hline
        + & 110001\\
        - & 111000\\
        - & 010101\\
        + & 101010\\
    \end{tabular}
    \end{minipage}
    }
    
    \medskip
    {\bf Case $n=10$.}
    
    \medskip
    \makebox[\textwidth][c]{
    \begin{minipage}[t]{0.25\textwidth}
    \centering
    \begin{tabular}{c|c}
        $\varepsilon_v$ & $v$ \\ \hline
        + & 0000011111\\
        - & 0000111110\\
        - & 0001111100\\
        + & 0011111000\\
        - & 0111110000\\
        + & 1000001111\\
        - & 1100000111\\
        + & 1110000011\\
    \end{tabular}
    \end{minipage}
    \begin{minipage}[t]{0.25\textwidth}
    \centering
    \begin{tabular}{c|c}
        $\varepsilon_v$ & $v$ \\ \hline
        + & 1111000001\\
        - & 1111100000\\
        + & 0001011101\\
        - & 0010111010\\
        - & 0100010111\\
        + & 0101110100\\
        + & 0111010001\\
        - & 1000101110\\
    \end{tabular}
    \end{minipage}
    \begin{minipage}[t]{0.25\textwidth}
    \centering
    \begin{tabular}{c|c}
        $\varepsilon_v$ & $v$ \\ \hline
        - & 1010001011\\
        + & 1011101000\\
        + & 1101000101\\
        - & 1110100010\\
        + & 0010011011\\
        - & 0011011001\\
        + & 0100110110\\
        - & 0110010011\\
    \end{tabular}
    \end{minipage}
    \begin{minipage}[t]{0.25\textwidth}
    \centering
    \begin{tabular}{c|c}
        $\varepsilon_v$ & $v$ \\ \hline
        + & 0110110010\\
        - & 1001001101\\
        + & 1001101100\\
        - & 1011001001\\
        + & 1100100110\\
        - & 1101100100\\
        - & 0101010101\\
        + & 1010101010\\
    \end{tabular}
    \end{minipage}
    }

    \medskip
    {\bf Case $n=12$.}
    
    \medskip
    \makebox[\textwidth][c]{
    \begin{minipage}[t]{0.25\textwidth}
    \centering
    \begin{tabular}{c|c}
        $\varepsilon_v$ & $v$ \\ \hline
+ & 000000111111\\
+ & 000001111110\\
+ & 000011111100\\
+ & 000111111000\\
+ & 001111110000\\
- & 011111100000\\
+ & 100000011111\\
- & 110000001111\\
- & 111000000111\\
- & 111100000011\\
- & 111110000001\\
- & 111111000000\\
- & 000010111101\\
- & 000101111010\\
+ & 001011110100\\
- & 010000101111\\
+ & 010111101000\\
+ & 011110100001\\
    \end{tabular}
    \end{minipage}
    \begin{minipage}[t]{0.25\textwidth}
    \centering
    \begin{tabular}{c|c}
        $\varepsilon_v$ & $v$ \\ \hline
- & 100001011110\\
- & 101000010111\\
+ & 101111010000\\
- & 110100001011\\
+ & 111010000101\\
+ & 111101000010\\
- & 000100111011\\
- & 001001110110\\
- & 001110110001\\
- & 010011101100\\
+ & 011000100111\\
+ & 011101100010\\
- & 100010011101\\
- & 100111011000\\
+ & 101100010011\\
+ & 110001001110\\
+ & 110110001001\\
+ & 111011000100\\
    \end{tabular}
    \end{minipage}
    \begin{minipage}[t]{0.25\textwidth}
    \centering
    \begin{tabular}{c|c}
        $\varepsilon_v$ & $v$ \\ \hline
- & 000110111001\\
- & 001000110111\\
- & 001101110010\\
+ & 010001101110\\
- & 011011100100\\
+ & 011100100011\\
- & 100011011100\\
+ & 100100011011\\
- & 101110010001\\
+ & 110010001101\\
+ & 110111001000\\
+ & 111001000110\\
+ & 000111000111\\
- & 001110001110\\
+ & 011100011100\\
- & 100011100011\\
+ & 110001110001\\
- & 111000111000\\
    \end{tabular}
    \end{minipage}
    \begin{minipage}[t]{0.25\textwidth}
    \centering
    \begin{tabular}{c|c}
        $\varepsilon_v$ & $v$ \\ \hline
+ & 001010110101\\
+ & 010010101101\\
+ & 010100101011\\
- & 010101101010\\
- & 010110101001\\
- & 011010100101\\
+ & 100101011010\\
+ & 101001010110\\
+ & 101010010101\\
- & 101011010100\\
- & 101101010010\\
- & 110101001010\\
+ & 001100110011\\
- & 011001100110\\
+ & 100110011001\\
- & 110011001100\\
- & 010101010101\\
+ & 101010101010\\
    \end{tabular}
    \end{minipage}
    }

\qed

We similarly deal with the case when $n$ is a small power of 2.
\begin{lem} \label{lem:npow2small}
    For $n\in \{2,4,8\}$, there exist signs $\varepsilon_v\in \{-1,1\}$ for each $v\in V_0$, and $w\in \ZZ^n$ with $w_i\geq -1$ for all $i$, such that $\sum_{v\in V_0} \varepsilon_v v=w$.
\end{lem}
{\bf Proof.}
    It is not hard to construct for $n=2$. For $n\in \{4,8\}$, as with the proof of Lemma~\ref{lem:nevensmall}, we use a SAT solver to assign signs to $v\in V_0^{\pm}$ whose orbits under $\rho$ also contains $-v$, such that $\sum \varepsilon_v v=2w$. 

    \medskip
    {\bf Case $n=4$.} Set $w=(3,-1,-1,-1)$ and the signs are given in the table below.
    
    \medskip
    \makebox[\textwidth][c]{
    \centering
    \begin{tabular}{c|c}
        $\varepsilon_v$ & $v$ \\ \hline
- & 0011\\
- & 0110\\
+ & 1001\\
+ & 1100\\
- & 0101\\
+ & 1010\\
    \end{tabular}
    }

    \medskip
    {\bf Case $n=4$.} Set $w=(3,-1,-1,-1,3,-1,-1,-1)$ and the signs are given in the table below.
    
    \medskip
    \makebox[\textwidth][c]{
    \begin{minipage}[t]{0.25\textwidth}
    \centering
    \begin{tabular}{c|c}
        $\varepsilon_v$ & $v$ \\ \hline
- & 00001111\\
+ & 00011110\\
- & 00111100\\
- & 01111000\\
+ & 10000111\\
+ & 11000011\\
- & 11100001\\
+ & 11110000\\
+ & 00101101\\
- & 01001011\\
+ & 01011010\\
    \end{tabular}
    \end{minipage}
    \begin{minipage}[t]{0.25\textwidth}
    \centering
    \begin{tabular}{c|c}
        $\varepsilon_v$ & $v$ \\ \hline
+ & 01101001\\
- & 10010110\\
- & 10100101\\
+ & 10110100\\
- & 11010010\\
- & 00110011\\
- & 01100110\\
+ & 10011001\\
+ & 11001100\\
- & 01010101\\
+ & 10101010\\
    \end{tabular}
    \end{minipage}
    }
    
\qed

\end{document}